\documentclass[12pt]{article}

\usepackage{graphicx}
\usepackage{enumerate}
\usepackage{amsmath}
\usepackage{amssymb}
\usepackage[latin1]{inputenc}
\def\r{\mathbb{R}}
\def\n{\mathbb{N}}
\def\c{\mathbb{C}}
\def\dd{\mathbb{D}}

\newtheorem{lemma}{Lemma}

\newtheorem{theorem}{Theorem}

\newtheorem{remark}{Remark}
\newtheorem{conjecture}{Conjecture}

\newcommand{\dist}{\operatorname{dist}}
\newcommand{\re}{\operatorname{Re}}
\newcommand{\im}{\operatorname{Im}}
\newcommand{\intc}{\operatorname{Int}}

\begin{document}

\title{On the existence of a proper minimal surface in $\r^3$ with the conformal type of a disk.\thanks{Research partially supported by DGICYT grant number BFM2001-3489.\newline 2000 Mathematics Subject Classification: Primary 53A10.\newline
Key words and phrases: Minimal Surfaces.}}
\author{Santiago Morales\vspace{2mm}}
\date{November 5, 2002}

\maketitle
\section{Introduction}

Proper minimal surfaces in $\r^3$ have peculiar properties that are not shared by general minimal surfaces; especially in the embedded case.

It has been proved that, under additional conditions, this family of surfaces has strong restrictions on their conformal structures. For instance, Huber and Osserman proved that if $M$ is a complete minimal surface with finite total curvature, then $M$ has the conformal type of a compact Riemann surface minus a finite number of points. In particular it is parabolic, that is, $M$ is not compact and $M$ does not carry a negative non-constant subharmonic function.

In the same context, Collin, Kusner, Meeks, and Rosenberg \cite{ckmr} have proved that if $M$ is a properly immersed  minimal surface in $\r^3$, then $M(+)=\{ (x_1,x_2,x_3) \in M \; : \; x_3 \geq 0 \}$ is parabolic.

These results have motivated the following conjecture: 
\begin{conjecture}[Meeks, Sullivan, \cite{meeks}]
If $f:M\rightarrow\r^3$ is a complete proper minimal immersion where $M$ is a Riemannian surface without boundary and with finite genus, then $M$ is parabolic.
\end{conjecture}

The main goal of this paper is to show a counterexample to the conjecture.

\begin{theorem}
There exists $\chi:\dd\longrightarrow\r^3$, a conformal proper minimal immersion defined on the unit disk.
\end{theorem}

The immersion $\chi$ is obtained as the limit of a sequence of minimal immersions with boundary. These boundaries go uniformly to infinity in such a way that we get properness at the limit.

The aforementioned sequence is constructed in a recursive way. In this process, the lemma in page \pageref{lemma} is crucial. In this lemma, we modify a given minimal surface $X$ near its boundary to obtain a new minimal surface $Y$ such that the norm of $Y$ along the boundary increases with respect to the norm of $X$ and, at same time, the norm of $Y$ is controlled by a large enough lower bound in a small neighborhood of the boundary.

The tools we have utilized in the proof of this lemma are those that Nadirashvili \cite{nadi} used to obtain a complete bounded minimal surface in $\r^3$: Runge's theorem and the López-Ros transformation.

It is important to remark that the geometry of the surface described in the theorem is very complicated; the convex hull of the image under $\chi$ of any closed set of $\dd$ containing an open arc of $\partial\dd$ is $\r^3$. 

I would like to point out that the same technique of this paper can be used to construct a complete minimal surface properly immersed in a ball of $\r^3$, \cite{morales}.

{\em Acknowledgments.} I would like to thank Francisco Martín for suggesting to me this line of work and for several informative conversations. I would also like to thank Francisco J. López for helpful criticisms of the paper.

\section{Background and Notation}
The purpose of this section is to fix the notation used in the paper and to summarize some results about minimal surfaces.

We set $D(z_0,r)=\{z\in\c\: :\: |z-z_0|<r\}$. By a {\em polygon} $P$ we mean a closed simple curve in $\r^2$ formed by a finite number of straight segments verifying $0\in\intc P$, where $\intc P$ denotes the interior domain bounded by the curve $P$. 

Let $X:D\rightarrow \r^3$ be a conformal minimal immersion defined on a simply connected domain $D$, and let $S=\{e_1,e_2,e_3\}$ be a set of orthogonal coordinates in $\r^3$. We label $(X(z))_{j,S}=\left<X(z),e_j\right>$, $j=1,2,3$. We write $(X(z))_j$ instead $(X(z))_{j,S}$ when it is clever which orthogonal frame we are using.

We define the Weierstrass representation of the minimal immersion $X$ in $S$, $\phi_{(X,S)}=\{\phi_{1,S},\phi_{2,S},\phi_{3,S}\}$, as
$$\phi_{j,S}(z)=\frac{\partial (X(z))_{j,S}}{\partial u} - i \frac{\partial (X(z))_{j,S}}{\partial v}, \qquad j=1,2,3, \quad z=u+iv,$$
The functions $\phi_{1,S},\phi_{2,S},\phi_{3,S}$ are holomorphic on $D$, verifying $\sum_{j=1}^3(\phi_{j,S})^2 \equiv 0$ and $\sum_{j=1}^3{| \phi_{j,S} |}^2 \not\equiv 0$. As usual, we define:  
$$f_{(X,S)}=\phi_{1,S}-i\phi_{2,S} \quad\text{and}\quad g_{(X,S)}=\frac{\phi_{3,S}}{\phi_{1,S}-i\phi_{2,S}}.$$
For $\phi$ the Weierstrass representation of $X$, we also denote $f_{(\phi,S)}=f_{(X,S)}$ and $g_{(\phi,S)}=g_{(X,S)}$.

Conversely, if we consider $S=\{e_1,e_2,e_3\}$ a set of Cartesian orthogonal coordinates in $\r^3$, and $f$ and $g$ are (respectively) a holomorphic and meromorphic functions on $D$ such that
$$\phi_{1,S}=\tfrac12f(1-g^2), \quad \phi_{2,S}=\tfrac i2f(1+g^2), \quad \phi_{3,S}=fg$$
are holomorphic functions on $D$, then
$$\begin{array}{c}
X:D \rightarrow \r^3,\\
\displaystyle X(z)= \sum_{j=1}^3 \left(\re \int_{z_0}^z \phi_{j,S}(w)dw\right)e_j+c, \qquad z_0 \in D,  \; c \in \r^3,
\end{array}$$
is a conformal minimal immersion.

We can write the conformal metric associated to the immersion $X$, $\lambda_X^2(z) \left< \cdot, \cdot \right>$, in terms of the Weierstrass representation as follows:
$$\lambda_X(z)=\tfrac12|f_{(X,S)}(z)|(1+|g_{(X,S)}(z)|^2)=\tfrac1{\sqrt2}\| \phi_{(X,S)} (z)\|.$$
Observe that the above formula does not depend on the orthogonal frame $S$.

\section{Proof of the Theorem}
In order to prove the theorem, we need the following lemma.

\begin{lemma}\label{lemma}
Let $X:O\longrightarrow\r^3$ be a conformal minimal immersion defined on a simply connected domain $O$ with $X(0)=0$. Consider $r>0$, $0<s<r/100$ and a polygon $P$ with $P\subset O$, satisfying:
\begin{equation}\label{eslabon}
r<\|X(z)\|<r+s/2,\quad \forall z\in O\setminus \intc P.
\end{equation}
Then, for any $b_1,b_2>0$, there exist a polygon $Q$ and a conformal minimal immersion $Y:U\longrightarrow \r^3$ defined on an open neighborhood of $\overline{\intc Q}$ with $Y(0)=0$, such that:
\begin{enumerate}[(a)]
\item $P\subset \intc Q\subset\overline{\intc Q}\subset U\subset O$;
\item $\|Y(z)-X(z)\|<b_1$, $\forall z\in\overline{\intc P}$;
\item $|\|Y(z)\|-(r+s)|<b_2$, $\forall z \in Q$;
\item $\|Y(z)\|> r-3\sqrt{sr}>r/2$, $\forall z \in \intc Q\setminus\intc P$.
\end{enumerate}
\end{lemma}

Roughly speaking, in this lemma, we modify a minimal surface $X$ near the boundary, (property {\em (b)}), to obtain a new minimal surface $Y$ such that the norm of $Y$ along a polygon increases with respect to the norm of $X$, (property {\em (c)}), and, at same time, the norm of $Y$ is controlled by a lower bound, (property {\em (d)}).  Properties {\em (c)} and {\em (d)} are crucial to obtain properness.

This lemma will be proved in Section \ref{sec:lemma}.

We use the lemma to construct two sequences $\{X_n\}_n$ and $\{P_n\}_n$, where $P_n$ is a polygon, and $X_n$ is a conformal minimal immersion defined on a neighborhood of $\overline{\intc P_n}$ with $X_n(0)=0$, satisfying the following properties for all $n\in\n$:
\begin{enumerate}[\bf ({T}1)$_{n}$]
\item $P_{n-1}\subset \intc P_n$ and $P_n \subset D(0,3)$;
\item $|\|X_n(z)\|-r_n|<\frac1{2(n+1)^2}$,  $\forall z\in P_n$;
\item $\|X_n(z)\|\geq \frac{r_{n-1}}2-\frac 1{2n}$, $\forall z\in \intc P_n\setminus\intc P_{n-1}$;
\item $\|X_n(z)-X_{n-1}(z)\|<\frac1{n^2}$, $\forall z\in\overline{\intc P_{n-1}}$;
\item $\lambda_{X_n}(z)\geq\alpha_n\lambda_{X_{n-1}}(z),$ $\forall z \in \overline{\intc P_{n-1}}$ where $\{\alpha_k\}_k$ is a sequence such that $0<\alpha_k<1$ and $\{\prod_{k=1}^j \alpha_k\}_j$ converges to $1/2$.
\end{enumerate}
where $r_k=r_{k-1}+2/k$ for $k>1$ and $r_1>301$.

Sequences $\{X_n\}$ and $\{P_n\}$ are constructed in a recursive way. We can take $X_1(u+iv)=r_1(u,v,0)$ and $P_1\subset D(0,3)$ a suitable polygon satisfying (T2)$_1$. Suppose that we have got $X_1,\ldots,X_{n-1}$ and $P_1,\ldots, P_{n-1}$.

Now we construct the $n^{\text{th}}$ term. We choose $\{\widehat\epsilon_k\}\searrow 0$, with $\widehat\epsilon_k<1/n^2$ for all $k$. For each $k$ we consider $Y_k:U_k\longrightarrow\r^3$ and $Q_k$, given by the lemma, for the following data: 
$$X=X_{n-1},\quad P=P_{n-1}, \quad r=r_{n-1}-\frac1n,\quad s=\frac3n,\quad b_1=\widehat\epsilon_k,\quad b_2=\frac1{2(n+1)^2},$$
and $O$ a simply connected domain with $\overline{\intc P_{n-1}}\subset O\subset D(0,3)$ and verifying (\ref{eslabon}). 
From {\em (b)} in the lemma, we deduce that the sequence $\{Y_k\}$ uniformly converges to $X_{n-1}$ on $\overline{\intc P_{n-1}}$. This implies that $\{\lambda_{Y_k}\}$ uniformly converges to $\lambda_{X_{n-1}}$  on $\overline{\intc P_{n-1}}$, and hence there is a $k_0 \in \n$ such that:
$$\lambda_{Y_{k_0}}(z) \geq \alpha_n \lambda_{X_{n-1}}(z), \qquad \forall z\in \overline{\intc P_{n-1}}.$$
We define $X_n=Y_{k_0}$, $P_n=Q_{k_0}$. It is easy to check (using the lemma) that $X_n$ and $P_n$ verify (T1)$_n$,\ldots, (T5)$_n$. This concludes the construction of the sequences $\{X_n\}$ and $\{P_n\}$.

Now, we define $\Delta=\bigcup_{n \in \n} \intc P_n$. $\Delta$ is a proper simply connected domain of $\c$, (see (T1)$_n$); thus $\Delta$ is biholomorphic to a disc. 

From (T4)$_n$, we have $\{X_n\}$ is a Cauchy sequence on compact sets in $\Delta$. Then there exists $\chi:\Delta\longrightarrow\r^3$ a harmonic map such that $\{X_n\}\rightarrow \chi$, uniformly on compact sets in $\Delta$. $\chi$ has the following properties:
\begin{itemize}
\item $\chi$ is minimal and conformal.
\item $\chi$ is an immersion. Indeed, for any $z \in \Delta$ there exists  $n_0 \in \n$ such that $z \in \intc P_{n_0}$. From Property (T5)$_n$ one has that,
$$ \lambda_{X_k}(z) \geq \alpha_k \lambda_{X_{k-1}}(z) \geq \ldots \geq \alpha_k \ldots \alpha_{n_0+1} \lambda_{X_{n_0}}(z) \geq \alpha_k \ldots \alpha_1 \lambda_{X_{n_0}(z)},\quad \forall k>n_0.$$
Taking the limit as $k \to \infty$, we deduce that
$$\lambda_\chi(z) \geq \tfrac12 \lambda_{X_{n_0}}(z) > 0,$$
and so $\chi$ is an immersion.
\item $\chi$ is proper. Let $B$ be a compact subset of $\r^3$. Consider $n_0$ such that for all $k>n_0$, one has 
\begin{equation}\label{morcilla}
\|p\|<\frac{r_{k-1}}2-\frac1{2k}-2, \quad \forall p\in B
\end{equation}
and, (using (T3)$_n$ and  (T4)$_n$), for all $z\in \intc P_k\setminus\intc P_{k-1}$:
\begin{equation}\label{chorizo}
\|\chi(z)\|\geq \frac{r_{k-1}}2-\frac1{2k}-\|X_k(z)-\chi(z)\|> \frac{r_{k-1}}2-\frac1{2k}-2.
\end{equation}
From (\ref{morcilla}) and (\ref{chorizo}), we have $\chi^{-1}(B)\cap(\intc P_k\setminus \intc P_{k-1})=\emptyset$, $\forall k>n_0$. Then $\chi^{-1}(B)\subset\intc P_{n_0}$. This implies that $\chi^{-1}(B)$ is compact.
\end{itemize}

This completes the proof of the theorem.

\section{Proof of the lemma}\label{sec:lemma}

To prove the lemma, we modify the immersion $X$ such that it changes almost anything on $\intc P$. The aim of this modification is increasing the norm of $X$ to get {\em(c)} and {\em (d)} in the lemma at the same time.

The proof of the lemma consists of two inductive process. In the first one, we obtain a new minimal immersion $X_n$ defined on a simply connected domain $\Omega$ containing $\overline{\intc P}$. The new immersion $X_n$ is close to $X$ on $\intc P$ and its norm is greater than $r+s$ at $a_1,\ldots,a_n$, a finite collection of points around the boundary of $\Omega$.

In the second process, from $X_n$ we obtain a new minimal immersion $Y$ defined on $\Omega$ that proves the lemma. To construct $Y$, we increase the norm of $X_n$ along $n$ curves of $\partial \Omega$ that join the points $a_i$. As a consequence, $Y$ verify {\em (c)} for $Q$ a polygon near $\partial\Omega$. Also, $Y$ is close to $X_n$ on $\intc P$, and then $Y$ and $X$ are close on $\intc P$.

In both processes, we need to control the immersions in such a way that {\em (d)} holds.
\vspace{1cm}

First, we prepare the first process. We are going to fix some constants and a finite collection of point $p_1,\ldots,p_n$ around the polygon $P$ that we will use along this proof. Our purpose is add poles at the points $p_i$ to the Weierstrass representation $\phi^0$ of $X$, in order to increase the norm of the immersion at the aforementioned points $a_1,\ldots,a_n$.

We define the following constants:
\begin{enumerate}[\em {1}.a\rm)]
\item $\lambda>1$ such that $\lambda^3<2$ and $\lambda(r+s)-2\lambda^3\sqrt{\lambda^4(r+s)^2-r^2}>r-2\sqrt{(r+s)^2-r^2}$;
\item $s'=+\sqrt{\lambda^4(r+s)^2-r^2}$;
\item $1>\epsilon_0>0$ chosen small enough so that certain inequalities that we are going to use in this paper, were true. $\epsilon_0$ only depends on $r,s,b_1,\lambda$ and $s'$.
\end{enumerate}

Let $P'$ be a new polygon and $W$ a simply connected open set such that  $P\subset \intc P'\subset\overline{\intc P'}\subset W\subset\overline W\subset O$.
For a small enough neighborhood $E$ of $P'$ (with $E\subset W\setminus\intc P$), we can define a continuous map $S$, where $S(z)$, $\forall z\in E$ is a set of orthogonal coordinates in $\r^3$, $\{e_1(z),e_2(z),e_3(z)\}$ with $e_3(z)=\frac{X(z)}{\|X(z)\|}$.\label{frame}

\begin{figure}[h] 
\begin{center}
\includegraphics[width=12.5cm]{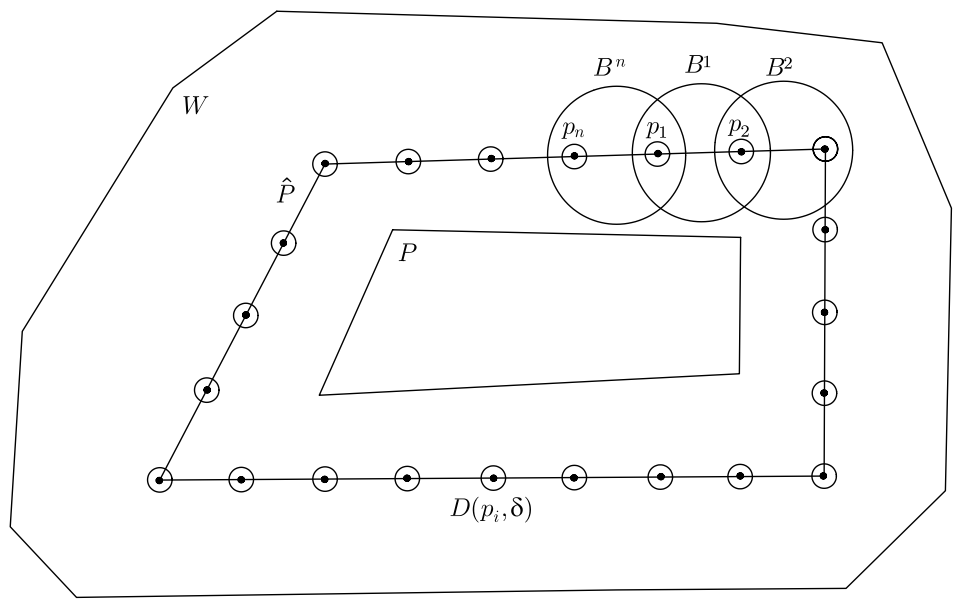}
\caption{Distribution of the points $p_i$.}\label{lospe}
\end{center}
\end{figure}

We choose a finite collection of point $\{p_1,\ldots,p_n\}\subset E$, (we label $p_{n+1}=p_1$), such that the segments $\overline{p_1p_2},\ldots,\overline{p_{n-1}p_n},\overline{p_np_{n+1}}$ forme a new polygon $\widehat P\subset E$ and verifying the following properties:
\begin{enumerate}[\em {2}.a\rm)]
\item for all $i=1,\ldots,n$, there exists a disc $B^i\subset E$ with $p_i,p_{i+1}\in B^i$ and small enough so that $X$ is close to a constant map on $B^i$, i.e.
\begin{equation}\label{tarariX}
\|X(z)-X(w)\|<\epsilon_0, \quad \forall z, w \in B^i;
\end{equation}
\item for all $i=1,\ldots,n$, the sets of orthogonal coordinates $S(p_i)$ and $S(p_{i+1})$ are close, i.e.
\begin{equation}\label{tarariS}
\left\|e_j^i-e_j^{i+1}\right\|<\epsilon_0,\quad \forall j\in\{1,2,3\};
\end{equation}
where we label $e^i_j=e_j(p_i)$.
\item for all $i=1,\ldots,n$,
\begin{equation}\label{tararif}
f_{(X,S(p_i))}(p_i)\not=0.
\end{equation}
We can choose $p_1,\ldots,p_n$ verifying (\ref{tararif}) because $f_{(X,S(p))}(p)=0$ can not hold on any open subset. If we have $f_{(X,S(p))}(p)=0$ on an open subset then, $X(p)$ is parallel to the normal vector of $X$ at $p$, $X(p)=\|X(p)\|N(p)$, but this is impossible on an open set for $X$ minimal immersion.
\end{enumerate}
Note that $n$ diverges as $\epsilon_0\rightarrow 0$.

We define $n$ complex numbers $\theta_1,\ldots,\theta_n$, with $|\theta_i|=1$ and $\im \theta_i\not=0$ such that
\begin{equation}\label{calambre}
\left| \overline{\theta_i}\frac{\overline{f_{(X,S(p_i))}(p_i)}}{|f_{(X,S(p_i))}(p_i)|}-\overline{\theta_{i+1}}\frac{\overline{f_{(X,S(p_{i+1}))}(p_{i+1})}}{|f_{(X,S(p_{i+1}))}(p_{i+1})|}\right|<\epsilon_0,\quad \forall i=1,\ldots,n.
\end{equation}
Remark that {\em 2.a}), {\em 2.b}) and (\ref{calambre}) holds for $i=n$ where $p_{n+1}=p_1$. 

For every point $p_i$, $i=1,\ldots,n$, we consider a disk $D(p_i,\delta)$ where $0<\delta<\epsilon_0^2$ chosen small enough so that, (see Figure \ref{lospe}):
\begin{enumerate}[\em {3}.a\rm )]
\item $\overline{\intc \widehat P\setminus \cup_{k=1}^n D(p_k,\delta)}$ is a simply connected set;
\item $\overline{D(p_i,\delta)\cup D(p_{i+1},\delta)}\subset B^i$, $\forall i=1,\ldots,n$;
\item $D(p_i,\delta)\cap D(p_k,\delta)=\emptyset$, $\forall i=1,\ldots,n$, and $k\not=i$;
\item $\sqrt\delta\max_{\overline{D(p_i,\delta)}}\{|f_{(X,S(p_i))}|\} <1$, $\forall i=1,\ldots,n$;
\item $\sqrt\delta\frac{\max_{\overline{D(p_i,\delta)}} \{|f_{(X,S(p_i))}g^2_{(X,S(p_i))}|\}}{|\im\theta_i|} <1$, $\forall i=1,\ldots,n$;
\item $\sqrt\delta\max_{\overline{D(p_i,\delta)}} \{\|\phi^0\|\}<1$, $\forall i=1,\ldots, n$.
\end{enumerate}
We finish these previous steps defining $l$ as
\begin{equation}\label{quefrio}
\begin{split}
l=&\sup_{z\in \overline{\intc \widehat P\setminus \cup_{k=1}^n D(p_k,\delta)}} \inf\{ \text{Euclidean length of } \alpha \: :\: \alpha \text{ is a curve in }\\
& \overline{\intc \widehat P\setminus \cup_{k=1}^n D(p_k,\delta)}
\text{ with origin 0 and ending at } z\}+2\pi\delta+\delta+1.
\end{split}
\end{equation}
Observe that $\delta$ depend on $\epsilon_0$ and $\theta_1,\ldots,\theta_n$.

\subsection{The first inductive process}
In this process, we will modify the Weierstrass date of $X$, $\phi^0$, to produce a sequence of $n$ new Weierstrass data $\phi^1,\ldots,\phi^n$. In the $i^{\text{th}}$ step, we modify $\phi^{i-1}$ on the disk $D(p_i,\delta)$ to obtain a new Weierstrass representation $\phi^i$ with a pole at $p_i$ that is close to $\phi^{i-1}$ outside some neighborhoods of $p_1,\ldots,p_i$, (see next Property (A6)$_i$). To do this we will use the López-Ros transformation in the orthogonal frame $S(p_i)$. At the same time, we have to have a control on $\phi^i$ along the segment $\overline{p_i (p_i+\delta)}$ (see property (A7)$_i$, (A8)$_i$ and (A9)$_i$). 

When we have $\phi^1,\ldots,\phi^n$, we are going to define a simply connected domain $\Omega$ that contains $\overline{\intc P}$ and does not contain any poles of $\phi^1,\ldots,\phi^n$. Then $X_i(z)=\int_0^z \phi^i$, $z\in\Omega$ is well defined for all $i=1,\ldots,n$. From (A6)$_i$, we have that $X_i$ and $X_{i-1}$ are close outside $D(p_i,\delta)$ (see property (P1)$_i$). In particular $X_n$ are close to $X_1$ outside $\Omega\setminus\cup_{k=1}^nD(p_k,\delta)$. As $\phi^i$ has a pole at $p_i$, we obtain a point $a_i$ that is close to $p_i$ and verifies (A7)$_i$, that is, $\|X_n(a_i)\|$ is greater than $r+s$, (see property (P4)$_i$). The control that we have on $\phi^i$ along the segment $\overline{p_i (p_i+\delta)}$ implies that $X_n(a_i)$ is close to $X_n(a_{i+1})$ (see property (P3)$_i$). This fact will be crucial along the second inductive process, Section \ref{seccion pancrisima}.
\vspace{.8cm}

\begin{figure}[hbtp] 
\begin{center}
\includegraphics[width=5.5cm]{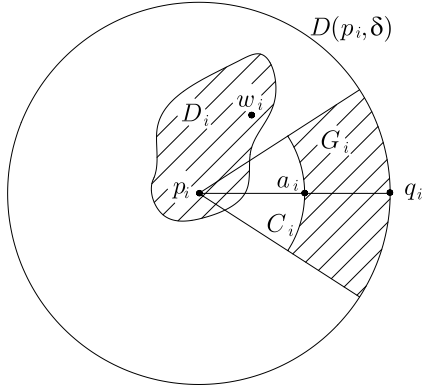}
\caption{The points $q_i$ and $a_i$.}\label{que}
\end{center}
\end{figure}

We are going to construct in a recursive way a sequence $\Psi_i=\{\phi^i,k_i,a_i,C_i,G_i,D_i\}$, $i=1,\ldots,n$, where:
\begin{enumerate}[\em {4}.a\rm )]
\item $\phi^i:\overline W\rightarrow \c^3$ is a Weierstrass representation with poles at $p_1,\ldots,p_i$;
\item $k_i$ is a suitable positive constant;
\item $a_i$ is a point lying on the segment $\overline{p_iq_i}$, where $q_i=p_i+\delta$;
\item $C_i$ is an open arc of a circumference centered at $p_i$ with $a_i\in C_i$;
\item $G_i$ is the closed annular sector bounded by $C_i$, a piece of $\partial D(p_i,\delta)$ and two radii as Figure \ref{que} indicates;
\item $D_i$ is an open simply connected subset of $\c$ verifying $\overline{D_i}\cap G_i=\emptyset$ and
$$\{p_i,w_i=p_i-k_i\theta_i\}\subset D_i\subset\overline{D_i}\subset D(p_i,\delta).$$
\end{enumerate}
\begin{remark}\label{n+1}
In what follows, we will use the convention that $\Psi_{n+1}\stackrel{\text{\tiny def}}= \Psi_1$.
\end{remark}
The above sequence is constructed in order to satisfy the following properties:

\begin{enumerate}[\bf ({A}1)$_{i}$]
\item $\sqrt\delta\max_{\overline{D(p_k,\delta)}}\{|f_{(\phi^i,S(p_k))}|\}<1$, $k=i+1,\ldots,n$;
\item $\sqrt \delta\frac{\max_{\overline{D(p_k,\delta)}} \{|f_{(\phi^i,S(p_k))}g^2_{(\phi^i,S(p_k))}|\}}{|\im\theta_k|}<1$, $k=i+1,\ldots,n$;
\item $2\lambda^2s'\frac{|f_{(\phi^i,S(p_k))}(p_k)-f_{(\phi^0,S(p_k))}(p_k)|}{|f_{(\phi^0,S(p_k))}(p_k)|}<\epsilon_0$, $k=i+1,\ldots,n$;
\item $\|\re\int_{\alpha_z}\phi^i\|<\epsilon_0$, $\forall z\in C_i$ where $\alpha_z$ is a piece of $C_i$ joining $a_i$ and $z$;
\item In the orthogonal frame $S(p_i)$, one has $\phi_3^i=\phi_3^{i-1}$;
\item $\|\phi^i(z)-\phi^{i-1}(z)\|<\frac{\epsilon_0}{n l}$, $\forall z\in\overline W\setminus (D(p_i,\delta)\cup(\cup_{k=1}^{i-1} D_k))$;
\item Using the frame $S(p_i)$, we have $\|(\re\int_{\overline{q_ia_i}} \phi^i_1)e_1^i+(\re\int_{\overline{q_ia_i}} \phi^i_2)e_2^i\|>\lambda s'$; 
\item In $S(p_i)$, one has $\|(\re\int_{\overline{q_i z}} \phi^i_1)e_1^i+(\re\int_{\overline{q_i z}} \phi^i_2)e_2^i\|<\lambda^3s'$, $\forall z\in G_i$;
\item $\|\re\int_{\overline{q_ia_i}}\phi^i-\re\int_{\overline{q_{i-1}a_{i-1}}}\phi^{i-1}\|<\epsilon_0(4\lambda^2s'+12)$, $i=2,\ldots,n+1$.
\end{enumerate}
All the above properties are true for $i=1,\ldots,n$, except for Properties (A1)$_i$, (A2)$_i$ and (A3)$_i$ that hold only for $i=1,\ldots,n-1$. In the same way, Property (A9)$_i$ is valid for $i=2,\ldots,n+1$ (see Remark \ref{n+1}).

As we have mentioned at the beginning of this section, we proceed by recursion.

Let $\phi^0$ be the Weierstrass representation of $X=X_0$.
Assume we have constructed $\Psi_1,\ldots,\Psi_{i-1}$. We are going to construct $\Psi_i$.

Define the meromorphic data $\phi^i$, in the orthogonal frame $S(p_i)$, as follows: 
$$f_{(\phi^i,S(p_i))}=f_{(\phi^{i-1},S(p_i))}h_i,\quad g_{(\phi^i,S(p_i))}=g_{(\phi^{i-1},S(p_i))}/h_i,$$
where
$$h_i(z)=\frac{k_i\theta_i}{z-p_i}+1.$$
Furthermore, we assume that $k_i>0$ is small enough to verify:
\begin{equation}\label{qui}
k_i\max_{w\in\overline{D(p_i,\delta)}}\left\{\left|\frac{f_{(\phi^{i-1},S(p_i))}(w)-f_{(\phi^{i-1},S(p_i))}(p_i)}{w-p_i}\right|\right\}<1,
\end{equation}
and (A1)$_i$, (A2)$_i$, (A3)$_i$ and (A6)$_i$. It is possible because
$\phi^i\rightarrow\phi^{i-1}$, as $k_i\rightarrow 0$, uniformly on $\overline W\setminus (D(p_i,\delta)\cup(\cup_{k=1}^{i-1} D_k))$ and (A1)$_{i-1}$, (A2)$_{i-1}$ and (A3)$_{i-1}$ hold. In the case $i=1$, notice that we have chosen $\delta$ such that {\em 3.d}) and {\em 3.e}) were true.

\begin{remark}
The meromorphic function $h_i$ is close to 1 outside a neighborhood of $p_i$. The constant $\theta_i$ has the effect of a rotation over $h_i$, when $z$ is close to $p_i$. Outside a neighborhood of $p_i$ this "rotation effect" almost disappears.
\end{remark}

From the definition of $\phi^i$, (A5)$_i$ trivially holds.

We define $a_i$ as the first point in $\overline{q_ip_i}$ when we move along this segment from $q_i$ to $p_i$, and such that:
\begin{equation}\label{IA}
\tfrac12|f_{(\phi^0,S(p_i))}(p_i)|\int_{\overline{q_ia_i}} \frac{k_idw}{w-p_i}=\lambda^2s'.
\end{equation}
Observe that $f_{(\phi^0,S(p_i))}(p_i)\not=0$ and $\int_{\overline{q_iz}}\frac{k_idw}{w-p_i}\in\r^+$, $\forall z\in\overline{q_ip_i}$.

Our next step consist of seeing the following:
\begin{multline}\label{trues}
\left\|\left(\re\int_{\overline{q_iz}}\phi_1^i(w)dw\right)e_1^i+\left(\re\int_{\overline{q_iz}}\phi_2^i(w)dw\right)e_2^i-\right.\\
\left.-\tfrac12\left(\int_{\overline{q_iz}}\frac{k_idw}{w-p_i}\right)\left(\re\overline{f_{(\phi^0,S(p_i))}(p_i)\theta_i}\:e_1^i+\im\overline{f_{(\phi^0,S(p_i))}(p_i)\theta_i}\:e_2^i\right)\right\|<2\epsilon_0,
\end{multline}
where $z\in\overline{q_ia_i}$, and $\phi^i_1$ and $\phi^i_2$ are expressed in the orthogonal frame $S(p_i)$.

In other words, the first two coordinates in $S(p_i)$ of the curve $z\mapsto \re\int_{\overline{q_iz}}\phi^i(w)dw$, $z\in\overline{q_ia_i}$, approximates the segment starting at 0 in the direction of $\re\overline{f_{(\phi^0,S(p_i))}(p_i)\theta_i}\: e_1^i+\im\overline{f_{(\phi^0,S(p_i))}(p_i)\theta_i}\: e_2^i$.

To get (\ref{trues}), for $z\in\overline{q_ia_i}$, we write:
\begin{equation}\nonumber
\begin{split}
&\left(\re\int_{\overline{q_iz}}\phi_1^i(w)dw\right)+i\left(\re\int_{\overline{q_iz}}\phi_2^i(w)dw\right)=\\
&=\tfrac12\left(\int_{\overline{q_iz}}\overline{f_{(\phi^{i-1},S(p_i))}(w)h_i(w)dw}-\int_{\overline{q_iz}}f_{(\phi^{i-1},S(p_i))}(w)g^2_{(\phi^{i-1},S(p_i))}(w)\frac{dw}{h_i(w)}\right)=
\end{split}
\end{equation}
\begin{eqnarray}
&=&\tfrac12\int_{\overline{q_iz}}\overline{f_{(\phi^0,S(p_i))}(p_i)\frac{k_i\theta_i}{w-p_i}dw}+\nonumber\\
&+&\tfrac12\int_{\overline{q_iz}}\overline{(f_{(\phi^{i-1},S(p_i))}(p_i)-f_{(\phi^0,S(p_i))}(p_i))\frac{k_i\theta_i}{w-p_i}dw}+\label{LaVit1}\\
&+&\tfrac12\int_{\overline{q_iz}}\overline{(f_{(\phi^{i-1},S(p_i))}(w)-f_{(\phi^{i-1},S(p_i))}(p_i))\frac{k_i\theta_i}{w-p_i}dw}+\label{LaVit2}\\
&+&\tfrac12\int_{\overline{q_iz}}\overline{f_{(\phi^{i-1},S(p_i))}(w)dw}-\label{LaVit3}\\
&-&\tfrac12\int_{\overline{q_iz}}f_{(\phi^{i-1},S(p_i))}(w)g^2_{(\phi^{i-1},S(p_i))}(w)\frac{dw}{h_i(w)}\label{LaVit4}
\end{eqnarray}

To obtain an upper bound in (\ref{LaVit1}) we use the definition of $a_i$ and Property (A3)$_{i-1}$ in the following sense:
\begin{multline}\label{dd10}
\left|\int_{\overline{q_iz}}\overline{(f_{(\phi^{i-1},S(p_i))}(p_i)-f_{(\phi^0,S(p_i))}(p_i))\frac{k_i\theta_i}{w-p_i}dw}\right|\leq\\
\leq2\lambda^2s'\frac{|f_{(\phi^{i-1},S(p_i))}(p_i)-f_{(\phi^0,S(p_i))}(p_i)|}{|f_{(\phi^0,S(p_i)}(p_i)|}<\epsilon_0.
\end{multline}

Now, we deal with (\ref{LaVit2}). Taking (\ref{qui}) into account one has:
\begin{multline}\label{dd11}
\left|\int_{\overline{q_iz}}\overline{(f_{(\phi^{i-1},S(p_i))}(w)-f_{(\phi^{i-1},S(p_i))}(p_i))\frac{k_i\theta_i}{w-p_i}dw}\right|\leq \\
\leq \delta k_i \max_{w\in\overline{D(p_i,\delta)}}\left\{\left|\frac{f_{(\phi^{i-1},S(p_i))}(w)-f_{(\phi^{i-1},S(p_i))}(p_i)}{w-p_i}\right|\right\}<\delta<\epsilon_0.
\end{multline}

To get the same upper bound for (\ref{LaVit3}), we use (A1)$_{i-1}$,
\begin{equation}\label{dd12}
\left|\int_{\overline{q_iz}}\overline{f_{(\phi^{i-1},S(p_i))}(w)dw}\right|\leq \delta \max_{\overline{D(p_i,\delta)}}\{|f_{(\phi^{i-1},S(p_i))}|\}<\sqrt\delta<\epsilon_0.
\end{equation}

In order to bound (\ref{LaVit4}), observe that $|h_i(w)|=\sqrt{(\frac{k_i}{w-p_i})^2+2\re\theta_i(\frac{k_i}{w-p_i})+1}>$ $>\sqrt{1-\re^2\theta_i}=|\im\theta_i|$, $\forall w\in\overline{q_ip_i}$. Therefore, (A2)$_{i-1}$ leads us to
\begin{multline}\label{dd13}
\left|\int_{\overline{q_iz}}f_{(\phi^{i-1},S(p_i))}(w)g^2_{(\phi^{i-1},S(p_i))}(w)\frac{dw}{h_i(w)}\right| \leq\\
\leq \delta \frac{\max_{\overline{D(p_i,\delta)}} \{|f_{(\phi^{i-1},S(p_i))}g^2_{(\phi^{i-1},S(p_i))}|\}}{|\im\theta_i|}<\sqrt\delta<\epsilon_0.
\end{multline}

Inequalities (\ref{dd10}), (\ref{dd11}), (\ref{dd12}) and (\ref{dd13}) give (\ref{trues}).

From (\ref{trues}) and the definition of $a_i$ we deduce (A7)$_i$, and 
\begin{equation}\label{casia7}
\left\|\left(\re\int_{\overline{q_iz}}\phi_1^i\right)e_1^i+\left(\re\int_{\overline{q_iz}}\phi_2^i\right)e_2^i\right\|<\lambda^2 s'+2\epsilon_0<\lambda^3 s',\quad \forall z\in\overline{q_ia_i}.
\end{equation}

At this point, we are able to define $D_i$, $C_i$ and $G_i$.

Let $D_i$ be a simply connected domain with $\{p_i,p_i-k_i\theta_i\}\subset D_i\subset \overline{D_i}\subset D(p_i,\delta)$ and $\overline{D_i}\cap\overline{q_ia_i}=\emptyset$ (this is possible because $\im\theta_i\not=0$, thus we have $p_i-k_i\theta_i\not\in \overline{p_iq_i}$). Since (\ref{casia7}), we can take $C_i$ an open arc of a circumference centered at $p_i$ with $a_i\in C_i$ and small enough to verify (A4)$_i$ and (A8)$_i$ on $G_i$, where $G_i$ is the closed annular sector bounded by $C_i$, a piece of $\partial D(p_i,\delta)$ and two radii, as it is indicated in Figure \ref{que}.

To check (A9)$_i$, we write:
\begin{multline}\label{torres}
\left\|\re\int_{\overline{q_ia_i}}\phi^i-\re\int_{\overline{q_{i-1}a_{i-1}}}\phi^{i-1}\right\|= \\
=\left\|\sum_{k=1}^3 \left(\re\int_{\overline{q_ia_i}}\phi_{k,S(p_i)}^i\right)e_k^i-\sum_{k=1}^3\left(\re\int_{\overline{q_{i-1}a_{i-1}}}\phi_{k,S(p_{i-1})}^{i-1}\right)e_k^{i-1}\right\|\leq\\
\leq \sum_{k=1}^3 \left\|\left(\re\int_{\overline{q_ia_i}}\phi_{k,S(p_i)}^i\right)e_k^i-\left(\re\int_{\overline{q_{i-1}a_{i-1}}}\phi_{k,S(p_{i-1})}^{i-1}\right)e_k^{i-1}\right\|.
\end{multline}
Next step consists of getting an upper bound for the three addends in (\ref{torres}):

From (\ref{trues}), we have
\begin{multline}\nonumber
\left\|\left(\re\int_{\overline{q_ia_i}}\phi_{1,S(p_i)}^i(w)dw\right)e_1^i-\left(\re\int_{\overline{q_{i-1}a_{i-1}}}\phi_{1,S(p_{i-1})}^{i-1}(w)dw\right)e_1^{i-1}\right\|\leq\\
\leq\left\|\left(\tfrac12\re\overline{f_{(\phi^0,S(p_i))}(p_i)\theta_i}\int_{\overline{q_ia_i}}\frac{k_idw}{w-p_i}\right)e_1^i-\right.\\
-\left.\left(\tfrac12\re \overline{f_{(\phi^0,S(p_{i-1}))}(p_{i-1})\theta_{i-1}}\int_{\overline{q_{i-1}a_{i-1}}}\frac{k_{i-1}dw}{w-p_{i-1}}\right)e_1^{i-1}\right\|+4\epsilon_0=
\end{multline}
using (\ref{IA}),
\begin{equation}\nonumber
\begin{split}
=&\left\|\lambda^2s'\frac{\re \overline{f_{(\phi^0,S(p_i))}(p_i)\theta_i}}{|f_{(\phi^0,S(p_i))}(p_i)|}e_1^i-\lambda^2s'\frac{\re \overline{f_{(\phi^0,S(p_{i-1}))}(p_{i-1})\theta_{i-1}}}{|f_{(\phi^0,S(p_{i-1}))}(p_{i-1})|}e_1^{i-1}\right\|+4\epsilon_0\leq\\
&\quad\leq \lambda^2s'\left\|\frac{\re \overline{f_{(\phi^0,S(p_i))}(p_i)\theta_i}}{|f_{(\phi^0,S(p_i))}(p_i)|}e_1^i-\frac{\re \overline{f_{(\phi^0,S(p_{i-1}))}(p_{i-1})\theta_{i-1}}}{|f_{(\phi^0,S(p_{i-1}))}(p_{i-1})|}e_1^i\right\|+\\
&\quad+\lambda^2s'\left\|\frac{\re\overline{ f_{(\phi^0,S(p_{i-1}))}(p_{i-1})\theta_{i-1}}}{|f_{(\phi^0,S(p_{i-1}))}(p_{i-1})|}e_1^i-\frac{\re \overline{f_{(\phi^0,S(p_{i-1}))}(p_{i-1})\theta_{i-1}}}{|f_{(\phi^0,S(p_{i-1}))}(p_{i-1})|}e_1^{i-1}\right\|+4\epsilon_0\leq\\
&\quad\leq \lambda^2s'\left|\frac{\overline{f_{(\phi^0,S(p_i))}(p_i)\theta_i}}{|f_{(\phi^0,S(p_i))}(p_i)|}-\frac{\overline{f_{(\phi^0,S(p_{i-1}))}(p_{i-1})\theta_{i-1}}}{|f_{(\phi^0,S(p_{i-1}))}(p_{i-1})|}\right|+\\
&\quad+\lambda^2s'\|e_1^i-e_1^{i-1}\|\left|\frac{\re \overline{f_{(\phi^0,S(p_{i-1}))}(p_{i-1})}}{|f_{(\phi^0,S(p_{i-1}))}(p_{i-1})|}\right|+4\epsilon_0\leq
\end{split}
\end{equation}
(\ref{calambre}) and (\ref{tarariS}) apply and conclude:
$$\leq 2\lambda^2s' \epsilon_0+4\epsilon_0=\epsilon_0(2\lambda^2s'+4).$$

In the same way, we obtain a bound for the second addends in (\ref{torres}),
\begin{multline}\nonumber
\left\|\left(\re\int_{\overline{q_ia_i}}\phi_{2,S(p_i)}^i\right)e_2^i-\left(\re\int_{\overline{q_{i-1}a_{i-1}}}\phi_{2,S(p_{i-1})}^{i-1}\right)e_2^{i-1}\right\|\leq\\
\leq \left\|\left(\tfrac12\im \overline{f_{(\phi^0,S(p_i))}(p_i)\theta_i}\int_{\overline{q_ia_i}}\frac{k_idw}{w-p_i}\right)e_2^i-\right.\\
-\left.\left(\tfrac12\im \overline{f_{(\phi^0,S(p_{i-1}))}(p_{i-1})\theta_{i-1}}\int_{\overline{q_{i-1}a_{i-1}}}\frac{k_{i-1}dw}{w-p_{i-1}}\right)e_2^{i-1}\right\|+4\epsilon_0\leq \epsilon_0(2\lambda^2s'+4).
\end{multline}

To estimate the third addend in (\ref{torres}), we use the Properties (A6)$_k$, $k=1,\ldots,i-1$, and {\em 3.f}),
\begin{equation}\nonumber
\begin{split}
\left\|\left(\re\int_{\overline{q_ia_i}}\phi_{3,S(p_i)}^i\right)\right. & \left.e_3^i- \left(\re\int_{\overline{q_{i-1}a_{i-1}}}\phi_{3,S(p_{i-1})}^{i-1}\right)e_3^{i-1}\right\|\leq\\
\leq &\left|\re\int_{\overline{q_ia_i}}\phi_{3,S(p_i)}^{i-1}\right|+\left|\re\int_{\overline{q_{i-1}a_{i-1}}}\phi_{3,S(p_{i-1})}^{i-2}\right| \leq \\
\leq &\delta\left(\max_{\overline{D(p_i,\delta)}}\{\|\phi^{i-1}\|\}+\max_{\overline{D(p_{i-1},\delta)}}\{\|\phi^{i-2}\|\}\right)\leq\\
\leq &\delta \left(
\max_{\overline{D(p_i,\delta)}}\{\|\phi^0\|\}+\max_{\overline{D(p_{i-1},\delta)}}\{\|\phi^0\|\}+2\epsilon_0\right)\leq 2\sqrt\delta+2\delta\epsilon_0<4\epsilon_0.
\end{split}
\end{equation}

The above bounds imply (A9)$_i$.

Note that case $i=n+1$:
$$\left\|\left(\re\int_{\overline{q_1a_1}}\phi_{1,S(p_1)}^1(w)dw\right)e_1^1-\left(\re\int_{\overline{q_na_n}}\phi_{1,S(p_n)}^n(w)dw\right)e_1^n\right\|\leq\epsilon_0(4\lambda^2s'+12),$$
can be proved in the same way.

Afterwards, we have just defined a sequence of Weierstrass representation $\phi^0,$ $\phi^1,\ldots,\phi^n$ verifying the above properties ({A}1)$_i$,$\ldots$,({A}9)$_i$.

Now, we are going to take a domain $\Omega$ where immersions $X_i(z)\stackrel{\text{\tiny def}}=\re\int_0^z \phi^i(w)dw$, $\forall i=1,\ldots,n$ are well defined and then become minimal immersions. Later, we will describe some properties of these immersion that we will use in the second process.

\begin{figure}[hbtp] 
\begin{center}
\includegraphics[width=10.5cm]{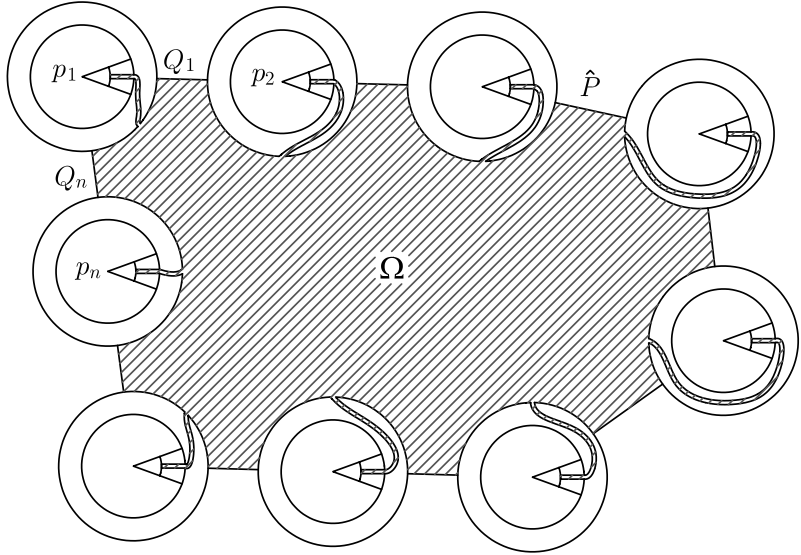}
\caption{The domain $\Omega$.}\label{ques}
\end{center}
\end{figure}

We define $\Omega$ as follows, (see figure \ref{ques}): Let $D^i$ be a disk centered at $p_i$ containing $D(p_i,\delta)$ for all $i=1,\ldots,n$, and let $\alpha_i\subset D^i\setminus \overline{D(p_i,\delta)}$ be a simple curve connecting $\partial D^i\cap \partial (\intc \widehat P\setminus \cup_{k=1}^n D^k)$ to the point $q_i$. Let $N_i$ be a small neighborhood of $\alpha_i\cup\overline{q_ia_i}$ in $\overline{G_i\cup (D^i\setminus D(p_i,\delta))}$. We define
$$\Omega=\left(\intc \widehat P\setminus\displaystyle\bigcup_{k=1}^n D^k\right)\cup\left(\displaystyle\bigcup_{k=1}^n N_k\right).$$

For suitable $D^i$, $\alpha_i$ and $N_i$, the following properties hold:
\begin{enumerate}[\em {5}.a\rm )]
\item $\overline \Omega$ is a simply connected domain;
\item $\overline{q_ia_i}\subset\overline\Omega$ and $\overline{\intc P}\subset\Omega$;
\item $\overline\Omega$ does not contain any $p_i$ point and any zero of $h_i$ for all $i=1,\ldots,n$;
\item $\sup_{z\in \overline{\Omega}} \inf \{$Euclidean length of $\alpha$ : $\alpha$ is a curve in $\overline\Omega$ with origin 0 and ending at $z\}<l$ where $l$ has been defined on (\ref{quefrio});
\item $\overline{\Omega}\cap \overline{D(p_i,\delta)}\subset G_i$.

\end{enumerate}
Taking {\em 5.a}) and {\em 5.c}) into account, we can define $n$ minimal immersions $X_1,\ldots,X_n$,
$$\begin{array}{c}
X_i:\Omega'\rightarrow\r^3\\ \\
\displaystyle X_i(z)=\re\int_0^z \phi^i(w)dw,\quad i=1,\ldots,n,
\end{array}$$
where $\Omega'$ is a suitable open neighborhood of $\overline \Omega$.

Consider an orthogonal frame $T_i=\{w_1^i,w_2^i,w_3^i\}$ for every $i=1,\ldots,n$ in such a way that
\begin{equation}\label{DonJavier}
w_3^i=\frac{X_n(a_i)}{\|X_n(a_i)\|}.
\end{equation}
These orthogonal frames will use in the second process.

We write $Q_i$, $i=1,\ldots,n$ as simple curves such that:
\begin{enumerate}[\em {6}.a\rm )]
\item $Q_i$ is connected, $\overline{Q_i}=Q_i$ and $Q_i\cap Q_j=\emptyset$ for $j\not=i$;
\item $Q_i$ connects $C_i$ with $C_{i+1}$;
\item $\partial \Omega$ is forme by the curves $Q_i$ and pieces of $C_i$;
\item $Q_i\subset B^i$;
\item $Q_i\cap \overline{D(p_k,\delta)}=\emptyset$ for $k\not\in\{i,i+1\}$;
\item We can suppose that $f_{(X_n,T_i)}(z)\not=0$, $\forall z\in Q_i$, $\forall i=1,\ldots,n$. To obtain this, we only need to make slight local modifications on the choice of $\Omega$. Note that the above modifications can be made in such a way that the properties of $\Omega$ and $Q_i$ remain.
\end{enumerate}

Now, we see some properties of the immersions $X_1,\ldots,X_n$ for $i=1,\ldots,n$:
\begin{enumerate}[\bf ({P}1)$_{i}$]
\item $\|X_i(z)-X_{i-1}(z)\|<\frac{\epsilon_0}n$, $\forall z\in\overline\Omega\setminus D(p_i,\delta)$;
\item $(X_i)_3=(X_{i-1})_3$, in the orthogonal frame $S(p_i)$;
\item $\|X_n(a_i)-X_n(a_{i+1})\|< \epsilon_0(4\lambda^2s'+17)$. Observe that $a_{n+1}=a_1$ when $i=n$;
\item $\|X_n(a_i)\|> \lambda(r+s)$;
\item $\|X_n(q_i)-X_n(z)\|<\lambda^3s'+5\epsilon_0$, $\forall z\in\overline{G_i\cap\Omega}$.
\end{enumerate}

Property (P1)$_i$ is a easy consequence of (A6)$_i$ and {\em 5.d}). On the other hand, it is straightforward to check (P2)$_i$ from (A5)$_i$.

Taking into account (P1), (A9)$_{i+1}$ and (\ref{tarariX}), we get (P3)$_i$ as follows:
\begin{equation}\nonumber
\begin{split}
\|X_n(a_i)-&X_n(a_{i+1})\|\leq \|X_i(a_i)-X_{i+1}(a_{i+1})\|+2\epsilon_0 \leq \\
&\leq \left\|\re\int_{\overline{q_ia_i}}\phi^i-\re\int_{\overline{q_{i+1}a_{i+1}}}\phi^{i+1}\right\|+\left\| X_i(q_i)-X_{i+1}(q_{i+1})\right\|+2\epsilon_0\leq\\
&\leq \epsilon_0(4\lambda^2s'+12)+ \| X_0(q_i)-X_0(q_{i+1})\|+4\epsilon_0\leq\\
&\leq\epsilon_0(4\lambda^2s'+12)+ 5\epsilon_0=\epsilon_0(4\lambda^2s'+17). 
\end{split}
\end{equation}

In order to prove (P4)$_i$, recall that
$$\|X_n(a_i)\|>\|X_i(a_i)\|-\epsilon_0.$$
We give bounds for the coordinates of $X_i(a_i)$ in $S(p_i)$ separately. Firstly, we bound the first and second coordinates.
\begin{multline}\nonumber
\|(X_i(a_i))_{1,S(p_i)}e_1^i+(X_i(a_i))_{2,S(p_i)}e_2^i\|>\\
>\left\|\left(\re\int_{\overline{q_ia_i}} \phi^i_1\right)e_1^i+\left(\re\int_{\overline{q_ia_i}} \phi^i_2\right )e_2^i\right\|
-\|(X_i(q_i))_1e_1^i+(X_i(q_i))_2e_2^i\|>
\end{multline}
using (A7)$_i$, (\ref{tarariX}) and the definition of $S(p_i)$ in page \pageref{frame},  we have
\begin{multline}\nonumber
>\lambda s'-\|(X_0(q_i))_1e_1^i+(X_0(q_i))_2e_2^i\|-\epsilon_0\geq\\
\geq\lambda s'-\|(X_0(p_i))_1e_1^i+(X_0(p_i))_2e_2^i\|-2\epsilon_0=\lambda s'-2\epsilon_0>s'.
\end{multline}
Now, we get a bound for the third coordinate. We use (P2)$_i$ and Hypothesis (\ref{eslabon}) in the lemma:
\begin{multline}\nonumber
|(X_i(a_i))_3|=|(X_{i-1}(a_i))_3|\geq |(X_0(a_i))_3|-\epsilon_0\geq\\
\geq |(X_0(p_i))_3|-2\epsilon_0= \|X_0(p_i)\|-2\epsilon_0> r-2\epsilon_0.
\end{multline}

Then, one has
\begin{multline}\nonumber
\|X_n(a_i)\|>\sqrt{(s')^2+( r-2\epsilon_0)^2}-\epsilon_0=\\
=\sqrt{\lambda^4(r+s)^2-4\epsilon_0(r-\epsilon_0)}-\epsilon_0>\lambda(r+s),
\end{multline}
that prove (P4)$_i$.

To obtain (P5)$_i$ we consider the orthogonal frame $S(p_i)$. Let $z\in\overline{G_i\cap\Omega}$. Using (A8)$_i$, we have 
\begin{multline}\nonumber
\|X_n(q_i)-X_n(z)\|\leq \|X_i(q_i)-X_i(z)\|+2\epsilon_0\leq\\
\leq\left\|\left(\re\int_{\overline{q_iz}}\phi_1^i\right)e_1^i+\left(\re\int_{\overline{q_iz}}\phi_2^i\right)e_2^i\right\|+|(X_i(q_i)-X_i(z))_{3,S(p_i)}|+2\epsilon_0<\\
< \lambda^3s'+|(X_0(q_i)-X_0(z))_{3,S(p_i)}|+4\epsilon_0\leq \lambda^3s'+5\epsilon_0.
\end{multline}

\subsection{The second inductive process}\label{seccion pancrisima}
In the first process, we have obtained a immersion $X_n$; and we have only guaranteed that its norm is greater than $r+s$ at the points $a_1,\ldots,a_n$, (Property (P4)). In this second process, our objective is increasing the norm of $X_n$ around the whole boundary of $\Omega$ in order to satisfy its norm is greater than $r+s$. So, we will modify the immersion $X_n$ to produce a sequence of minimal immersions $Y_1,\ldots,Y_n$ defined on $\overline\Omega$. In each step, we take certain small neighborhoods, $Q_i^\xi$ of the curve $Q_i$; and modify the Weierstrass representation of $Y_{i-1}$, $\psi^{i-1}$, on $Q_i^\xi$, to get a new minimal immersion $Y_i$. $Y_i$ and $Y_{i-1}$ are close on $\overline\Omega\setminus Q_i^\xi$, and the intrinsic metric of $Y_i$ increase on $Q_i^\xi$ with respect to $Y_{i-1}$. From this, we obtain that the norm of $Y_i$ is greater than $r+s$ around $Q_i$. To do this, we use the López-Ros transformation in the set of Cartesian coordinates $T_i$ and Runge's theorem.
So, note that the third coordinates of $X_n(a_i)$ and $X_n(a_{i+1})$, in the frame $T_i$,  are close and greater than $r+s$. This fact will be decisive in Section \ref{ultima}.

\begin{figure}[hbtp] 
\begin{center}
\includegraphics[width=14cm]{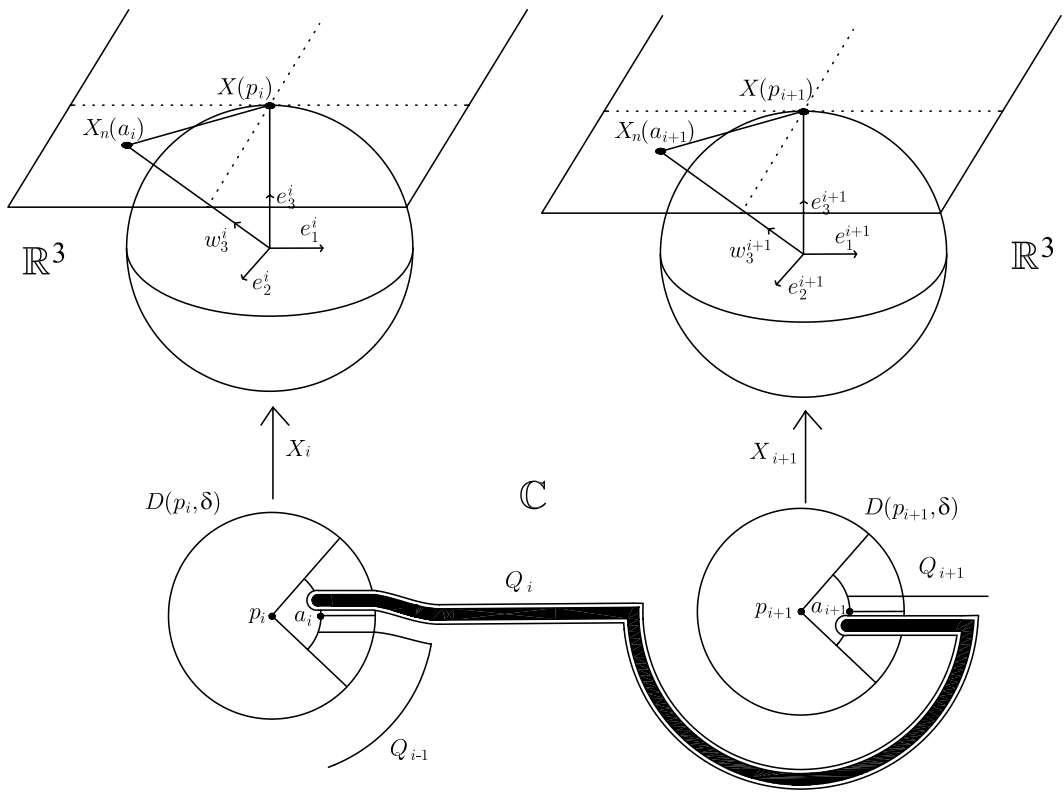}
\caption{Properties of $X_n$.}\label{idea2}
\end{center}
\end{figure}

Previously, we define a collection of neighborhoods of the arcs $C_i$ and a collection of neighborhoods of the curves $Q_i$:

For all $i=1\ldots,n$, let $\widehat C_i$ be an open neighborhood of $C_i$, small enough to satisfy:
\begin{equation}\label{mamanela}
\|X_n(z)-X_n(a_i)\|<3\epsilon_0,\quad \forall z\in \widehat C_i\cap\overline\Omega
\end{equation}
We can take $\widehat C_i$ verifying the above because, using (A4)$_i$, for all $z\in C_i$, $\|X_n(z)-X_n(a_i)\|< \|X_i(z)-X_i(a_i)\|+2\epsilon_0<3\epsilon_0$.

We define $0<\epsilon_1$ such that
\begin{equation}\label{calabera}
2\epsilon_1<\tfrac12\min_{Q_i}\{|f_{(X_n,T_i)}|\},\quad \forall i=1,\ldots,n.
\end{equation}
Observe that from {\em 6.f}), $\min_{Q_i}\{|f_{(X_n,T_i)}|\}\not=0$.

For all $i=1,\ldots,n$, we define $Q_i^{\xi}=\{z\in\c \: :\: \dist(z,Q_i)\leq\xi\}$, where 
$\dist(z,A)$ means the Euclidean distance between point $z$ and set $A$. We assume that $0<\xi$ is small enough to satisfy: 
\begin{enumerate}[\em {7}.a\rm)]
\item $Q_i^\xi \subset \Omega'$;
\item $Q_i^\xi\cap Q_j^\xi=\emptyset$, $i\not=j$;
\item $Q_i^\xi\cap\overline{D(p_k,\delta)}=\emptyset$, for $k\not\in\{i,i+1\}$;
\item $Q_i^{\xi}\subset B^i$;
\item $Q_i^{\xi/2}$ and $\overline{\Omega\setminus Q_i^\xi}$ are simply connected sets;
\item For all $z \in Q_i$ we have $|f_{(X_n,T_i)}(z)-f_{(X_n,T_i)}(x)|<\epsilon_1$, $\forall x\in B(z,\xi/2)$;
\item $\sup_{z\in\overline{\Omega\setminus Q_i^\xi}} \inf \{$Euclidean length of $\alpha$ : $\alpha$ is a curve in $\overline{\Omega\setminus Q_i^\xi}$ with origin 0 and ending at $z\}<l$.
\end{enumerate}

Label $Y_0\stackrel{\text{\tiny def}}=X_n$. We are going to construct in a recursive way $Y_1,\ldots,Y_n$, where $Y_i:\Omega'\rightarrow\r^3$ is a conformal minimal immersion with $Y_i(0)=0$ verifying the following properties:
\begin{enumerate}[\bf ({B}1)$_{i}$]
\item $(Y_i)_{3,T_i}=(Y_{i-1})_{3,T_i}$;
\item $\|Y_i(z)-Y_{i-1}(z)\|<\frac{\epsilon_0}n$, $\forall z\in\overline{\Omega\setminus Q_i^{\xi}}$;
\item $|f_{(Y_i,T_k)}(z)-f_{(Y_{i-1},T_k)}(z)|<\frac{\epsilon_1}n$, $\forall z\in \overline{\Omega\setminus Q_i^{\xi}}$, for $k=i+1
,\ldots,n$;
\item $\frac{\xi}2\left((\frac{1}{\tau_i}+\frac{\nu_i}{\tau_i(\tau_i-\nu_i)})\max_{Q_i^{\xi}}\{|f_{(Y_{i-1},T_i)}g^2_{(Y_{i-1},T_i)}|\}+\nu_i\max_{Q_i^{\xi}}\{|f_{(Y_{i-1},T_i)}|\}\right) <1$;
\item $\frac12\left( \frac{\tau_i\xi}{4}\min_{Q_i}\{|f_{(Y_0,T_i)}|\}-1\right)>2(r+s)+1$;
\end{enumerate}
where $\tau_i$ and $\nu_i$ are certain positive constants.

Assume we have $Y_0,Y_1,\ldots,Y_{i-1}$. We define $Y_i$ as:
$$f_{(Y_i,T_i)}=f_{(Y_{i-1},T_i)}l_i,\quad g_{(Y_i,T_i)}=g_{(Y_{i-1},T_i)}/l_i,$$
where $l_i:\c\rightarrow\c$ is a holomophic non zero function verifying:
\begin{itemize}
\item $|l_i(z)-\tau_i|<\nu_i$, $\forall z\in Q_i^{\xi/2}$;
\item $|l_i(z)-1|<\nu_i$, $\forall z\in\overline{\Omega\setminus Q_i^{\xi}}$.
\end{itemize}
The existence of this function \cite{nadi} is a consequence of Runge's theorem. We define $Y_i$ as the minimal immersion that becomes from the above meromorphic data.

We have assumed that $\tau_i>0$ is large enough and $\nu_i>0$ is small enough to satisfy (B2)$_i$, (B3)$_i$, (B4)$_i$ and (B5)$_i$. It is possible because the Weierstrass representation of $Y_{i-1}$ uniformly converges to the Weierstrass representation of $Y_i$, as $\nu_i\rightarrow 0$ and $\tau_i\rightarrow \infty$ on $\overline{\Omega\setminus Q_i^\xi}$. We take into account Property {\em 7.g}), to obtain (B2)$_i$.

Observe that (B1)$_i$ trivially holds.

\subsection{Definition of $Y$}\label{ultima}
Eventually, we define the open set $U$ as $U=\Omega$ and $Y:U\rightarrow\r^3$ as $Y=Y_n$. In this section, we are going to check that the immersion $Y$ verifies all the claims of the lemma.

Bearing in mind $\overline{\intc P}\subset \Omega\setminus ((\cup_{k=1}^n D(p_k,\delta))\cup(\cup_{k=1}^n Q_k^\xi))$, Property {\em (b)} is a easy consequence of (P1) and (B2).

The proof of the rest of claims of Lemma are nontrivial:

\paragraph{Properties {\em (a)} and {\em (c)}:} To prove the existence of a polygon $Q$ verifying {\em (a)} and {\em (c)}, we only need to get the following:
\begin{equation}\label{risitas}
\begin{split}
&\text{\em for all curve $\beta$ with origin 0 and ending at a}\\ 
&\text{\em point of $\partial\Omega$, there exists $z'\in\beta$ such that $\|Y(z')\|>r+s$.}
\end{split}
\end{equation}

We are going to see the above property. Let $\beta\subset\overline\Omega$ be a curve with $\beta(0)=0$ and $\beta(1)=z_0\in\partial\Omega$.

First, we want to show that $\|Y_n(z)\|>r+s$, $\forall z\in \widehat C_i$, for $i=1,\ldots,n$. Let $z\in\widehat C_i$, we study three possibilities: 
\begin{itemize}
\item Suppose that $z\in \widehat C_i\cap Q_i^\xi$, $i\in\{1,\ldots,n\}$. Then, considering $T_i$, and using (B2), (B1)$_i$, (\ref{mamanela}), (\ref{DonJavier}) and (P4)$_i$, we have
\begin{multline}\nonumber
\|Y_n(z)\|\geq\|Y_i(z)\|-\epsilon_0\geq|(Y_i(z))_3|-\epsilon_0=|(Y_{i-1}(z))_3|-\epsilon_0\geq |(Y_0(z))_3|-2\epsilon_0>\\
>|(Y_0(a_i))_3|-5\epsilon_0=\|Y_0(a_i)\|-5\epsilon_0>\lambda(r+s)-5\epsilon_0>r+s.
\end{multline}
\item Suppose that $z\in \widehat C_i\cap Q_{i-1}^\xi$. Now, considering $T_{i-1}$, one has 
$$\|Y_n(z)\|\geq |(Y_0(z))_{3,T_{i-1}}|-2\epsilon_0\geq$$
and from (\ref{mamanela}), (P3)$_{i-1}$ and (P4)$_{i-1}$, we deduce
\begin{multline}\nonumber
\geq |(Y_0(a_i))_{3,T_{i-1}}|-5\epsilon_0\geq  |(Y_0(a_{i-1}))_{3,T_{i-1}}|-\epsilon_0(4\lambda^2s'+17)-5\epsilon_0\geq\\
\geq \lambda(r+s)-\epsilon_0(4\lambda^2s'+17)-5\epsilon_0\geq r+s.
\end{multline}
Note that case $\widehat C_1\cap Q_n^\xi$ can be proved in the same way.
\item Suppose that $z\in \widehat C_i\setminus \cup_{k=1}^n Q_k^{\xi}$. Therefore,
$$\|Y_n(z)\|\geq \|Y_0(z)\| -\epsilon_0>\|Y_0(a_i)\|-4\epsilon_0>\lambda(r+s)-4\epsilon_0>r+s.$$
\end{itemize}

This proves that $\|Y_n(z)\|>r+s$, $\forall z\in\widehat C_i$. Hence, we can suppose that $\beta\cap\widehat C_i=\emptyset$, $\forall i=1,\ldots,n$.

From {\em 6.c}), $z_0\in Q_i$, $i=1,\ldots,n$. Let $z_1\in\beta\cap\partial B(z_0,\xi/2)$. Observe that $\overline{z_1z_0}\subset Q_i^{\xi/2}$. Next, we give a lower bound for $\|Y_n(z_0)-Y_n(z_1)\|$. We consider $T_i$, and write $f^{i-1}=f_{(Y_{i-1},T_i)}$, $g^{i-1}=g_{(Y_{i-1},T_i)}$. Then,
\begin{equation}\nonumber
\begin{split}
\|Y_n(z_0)-Y_n(z_1)\|&\geq \|Y_i(z_0)-Y_i(z_1)\|-2\epsilon_0\geq\\
&\geq \left|\left(\re\tfrac12\int_{\overline{z_1z_0}}f^{i-1}l_i-\frac{f^{i-1}(g^{i-1})^2}{l_i}\right)+\right.\\
&\left.+i\left(\re\tfrac i2\int_{\overline{z_1z_0}}f^{i-1}l_i+\frac{f^{i-1}(g^{i-1})^2}{l_i}\right)\right|-2\epsilon_0=\\
&=\tfrac12\left|\int_{\overline{z_1z_0}} \overline{f^{i-1}l_i}-\frac{f^{i-1}(g^{i-1})^2}{l_i}\right|-2\epsilon_0\geq\\
&\geq \tfrac{\tau_i}2\left|\int_{\overline{z_1z_0}} \overline{f^{i-1}}\right|-\tfrac1{2\tau_i}\left|\int_{\overline{z_1z_0}} f^{i-1}(g^{i-1})^2\right|-\\
&-\tfrac12\left|\int_{\overline{z_1z_0}} \overline{f^{i-1}(l_i-\tau_i)}\right|-\tfrac12\left|\int_{\overline{z_1z_0}} f^{i-1}(g^{i-1})^2\left(\frac1{l_i}-\frac1{\tau_i}\right)\right|-2\epsilon_0\geq\\
&\geq\tfrac{\tau_i}2\left|\int_{\overline{z_1z_0}} \overline{f^{i-1}}\right|-\tfrac\xi4\left(\tfrac{1}{\tau_i}\max_{Q_i^{\xi}}\{| f^{i-1}(g^{i-1})^2|\}+\right.\\
&\left.+\nu_i\max_{Q_i^{\xi}}\{|f^{i-1}|\}+\frac{\nu_i}{\tau_i(\tau_i-\nu_i)}\max_{Q_i^{\xi}}\{|f^{i-1}(g^{i-1})^2|\}\right)-2\epsilon_0\geq
\end{split}
\end{equation}
using (B4)$_i$, we continue
$$\geq \tfrac12\left(\tau_i\left|\int_{\overline{z_1z_0}} \overline{f^{i-1}}\right|-1\right)-2\epsilon_0.$$
On the other hand, from (B3)$_1,\ldots,$(B3)$_{i-1}$, {\em 7.f}), and (\ref{calabera}), we deduce
\begin{equation}\nonumber
\begin{split}
\left|\int_{\overline{z_1z_0}} \overline{f^{i-1}(w)dw}\right|&\geq\\
&\geq\left|f_{(Y_0,T_i)}(z_0)\int_{\overline{z_1z_0}} dw\right|-\left|\int_{\overline{z_1z_0}} (f_{(Y_0,T_i)}(z_0)-f_{(Y_0,T_i)}(w))dw\right|-\\
&-\left|\int_{\overline{z_1z_0}} (f_{(Y_0,T_i)}(w)-f^{i-1}(w))dw\right|\geq\\
&\geq \tfrac{\xi}2(|f_{(Y_0,T_i)}(z_0)|-\epsilon_1-\epsilon_1)\geq \tfrac{\xi}2(\min_{Q_i}\{|f_{(Y_0,T_i)}|\}-2\epsilon_1)\geq\\
&\geq \tfrac{\xi}4\min_{Q_i}\{|f_{(Y_0,T_i)}|\}.
\end{split}
\end{equation}
Therefore,
$$\|Y_n(z_0)-Y_n(z_1)\|\geq \tfrac12 \left( \tau_i\tfrac{\xi}{4}\min_{Q_i}\{|f_{(Y_0,T_i)}|\}-1\right)-2\epsilon_0\geq$$
and using (B5)$_i$, we obtain
\begin{equation}\label{extremadura}
\|Y_n(z_0)-Y_n(z_1)\|\geq 2(r+s)+1-2\epsilon_0>2(r+s).
\end{equation}
Indeed, one has $\|Y_n(z_0)\|>r+s$ or $\|Y_n(z_1)\|>r+s$.
 
This conclude the proof of (\ref{risitas}). Thus, there is a polygon $Q$ verify {\em (a)} and {\em (c)}.

\paragraph{Property {\em (d)}:} Let $z\in \intc Q\setminus\intc P$. We distinguish the five possible cases to prove this property:

\begin{itemize}
\item Suppose that $z\not\in (\cup_{k=1}^n D(p_i,\delta))\cup(\cup_{k=1}^n Q_i^{\xi})$. (P1), (B2) and (\ref{eslabon}) imply $\|Y_n(z)\|\geq \|X(z)\|-2\epsilon_0>r-2\epsilon_0>r-3\sqrt{sr}$.
\item Suppose that $z\in D(p_i,\delta)\setminus \cup_{k=1}^n Q_k^{\xi}$. We consider $S(p_i)$ and use (P1), (B2), (P2)$_i$, then:
\begin{multline}\nonumber
\|Y_n(z)\|\geq \|X_i(z)\|-2\epsilon_0\geq|(X_i(z))_3|-2\epsilon_0=\\
=|(X_{i-1}(z))_3|-2\epsilon_0\geq|(X(z))_3|-3\epsilon_0\geq
\end{multline}
and, from (\ref{tarariX}) and the definition of $S(p_i)$, one has
$$\geq |(X(p_i))_3|-4\epsilon_0=\|X(p_i)\|-4\epsilon_0>r-4\epsilon_0>r-3\sqrt{sr}.$$

\item Suppose that $z\in D(p_i,\delta)\cap Q_i^{\xi}$. Now, we consider the orthogonal frame $T_i$, and obtain
\begin{multline}\nonumber
\|Y_n(z)\|\geq\|Y_i(z)\|-\epsilon_0\geq|(Y_i(z))_3|-\epsilon_0=|(Y_{i-1}(z))_3|-\epsilon_0 \geq |(Y_0(z))_3|-2\epsilon_0\geq\\
\geq|(Y_0(a_i))_3|-|(Y_0(a_i))_3-(Y_0(q_i))_3|-|(Y_0(q_i))_3-(Y_0(z))_3|-2\epsilon_0>
\end{multline}
Note that $|(Y_0(a_i))_{3,T_i}|=\|Y_0(a_i)\|$. Since $z\in G_i$, we can use (P4)$_i$ and (P5)$_i$ to get the following
$$>\lambda(r+s)-2\lambda^3s'-12\epsilon_0=\lambda(r+s)-2\lambda^3\sqrt{\lambda^4(r+s)^2-r^2}-12\epsilon_0>$$
next we recall the definition of $\lambda$ on {\em 1.b}), and $s<r/100$, then
$$>r-2\sqrt{(r+s)^2-r^2}-12\epsilon_0>r-3\sqrt{sr}.$$

\item Suppose that $z\in D(p_{i+1},\delta)\cap Q_i^{\xi}$. We consider $T_i$. This case is similar to the former case:
\begin{multline}\nonumber
\|Y_n(z)\|\geq|(Y_0(a_i))_{3,T_i}|-|(Y_0(a_i))_{3,T_i}-(Y_0(a_{i+1}))_{3,T_i}|-\\
-|(Y_0(a_{i+1}))_{3,T_i}-(Y_0(q_{i+1}))_{3,T_i}|-|(Y_0(q_{i+1}))_{3,T_i}-(Y_0(z))_{3,T_i}|-2\epsilon_0\geq
\end{multline}
Here, (P3)$_i$ gives
$$>\lambda(r+s)-\epsilon_0(4\lambda^2+17)-2\lambda^3s'-12\epsilon_0\geq r-3\sqrt{sr}.$$
Note that case $z\in D(p_1,\delta)\cap Q_n^\xi$ can be proved in the same way, considering $T_n$.

\item Suppose that $z\in Q_i^{\xi}\setminus\cup_{k=1}^n D(p_k,\delta)$.
We consider $T_i$. A similar computation gives
\begin{multline}\nonumber
\|Y_n(z)\|\geq |(Y_0(a_i))_3|-|(Y_0(a_i))_3-(Y_0(q_i))_3|-|(Y_0(q_i))_3-(Y_0(z))_3|-2\epsilon_0\geq\\
\geq\lambda(r+s)-\lambda^3s'-7\epsilon_0-|(Y_0(q_i))_3-(Y_0(z))_3|
\end{multline}
and using properties (P1) and (\ref{tarariX}), we obtain
$$\geq\lambda(r+s)-\lambda^3s'-10\epsilon_0\geq r-3\sqrt{sr}.$$
\end{itemize} 

This completes the proof of the lemma.

\noindent
{\itshape
Santiago Morales \\
Departamento de Geometría y Topología \\
Universidad de Granada \\
18071 Granada, Spain \\
E-mail: {\upshape\ttfamily santimo@ugr.es} \\
URL: {\upshape\ttfamily http://www.ugr.es/\~{ }santimo }   }


\begin{thebibliography}{99}
\bibitem{ckmr} P. Collin, R. Kusner, W. H. Meeks III, and H. Rosenberg {\em The topology, geometry and conformal structures of properly embedded minimal surfaces}. In preparation.

\bibitem{morales} F. Martín and S. Morales, {\em A note about the asymptotic behavior of a complete bounded minimal surface in $\r^3$}. Preprint.

\bibitem{meeks} W. H. Meeks III, {\em Global problems in Classical minimal surfaces theory}. In preparation.



\bibitem{nadi}N. Nadirashvili, {\em Hadamard's and Calabi-Yau's conjectures on negatively curved and minimal surfaces.}
Invent. math., 126 (1996), 457-465.
\end{thebibliography}
\end{document}